\documentclass[11pt]{article}
\usepackage{graphicx}
\usepackage{amssymb,amsmath}
\usepackage[isolatin]{inputenc}
\usepackage{pdfsync}
\newtheorem{teorema}{Theorem}[section]
\newtheorem{prop}[teorema]{Proposition}
\newtheorem{coro}[teorema]{Corollary}
\newtheorem{lema}[teorema]{Lemma}
\newtheorem{defi}{Definition}[section]

\newcommand{\snn}{\operatorname{sn}}
\newcommand{\sn}{\snn_c}
\newcommand{\cnn}{\operatorname{cn}}
\newcommand{\cn}{\cnn_c}
\newcommand{\ctt}{\operatorname{cot}}
\newcommand{\ctn}{\ctt_c}

\newtheorem{remark}[teorema]{Remark}

\usepackage{makeidx,multicol}

\makeindex 

\def\Z{\mathbb Z}

\def\R{\mathbb R}

\newcommand{\hh}{\mathbb{H}}
\newcommand{\Ss}{\mathbb{S}}

\title{Evolutes and isoperimetric deficit in two-dimensional spaces of constant curvature}
\author{Julià Cufí and Agustí Reventós}
\date{}
\begin{document} \maketitle 



\abstract{We relate the total curvature and the isoperimetric deficit of a curve $\gamma$ in a two-dimensional space of constant curvature with the area enclosed by the evolute of $\gamma$.  We provide also a Gauss-Bonnet theorem for a special class of evolutes.\footnote{Work partially supported by grants MTM2012-36378 and MTM2012-34834 (MEC).}

}
\section{Introduction} The setting of this paper is the space  
$X^2_c$,  the $2$-dimensional complete and simply connected riemannian  manifold   of constant curvature $c$, i.e.  the sphere  $ \Ss^2_c$ of radius $R=\frac{1}{\sqrt{c}}$ for $c>0$,  the hyperbolic plane $\hh^2_c$ for $c<0$ (the imaginary sphere of radius $Ri=\frac{1}{\sqrt{c}}$), or the Euclidean plane for $c=0$. We shall assume $X^2_c$ oriented. 

For a closed curve $\gamma$ on $X_{c}^2$ we will consider the evolute $\gamma_{e}$ of $\gamma$ and denote by $F_{e}$ the area with multiplicities  enclosed by $\gamma_{e}$. By means of $|F_{e}|$ we estimate the deficit of the total curvature and the isoperimetric deficit of the curve $\gamma$.

The integral of the curvature of a simple closed curve (the total curvature) in the Euclidean space $\R^3$ has been widely studied. The most remarkable result is Fenchel's  theorem which states that this integral is greater than, or equal to, $2\pi$. It is equal to $2\pi$ if and only if the curve is a plane convex curve; see \cite{Fenchel}. The following result gives an interpretation of the difference between the total curvature and $2\pi$ for curves on $X_{c}^2$.

\begin{teorema}
  Let $\gamma(s)$ be a positively oriented closed strongly convex  curve on $X^2_c$ parametrized by arclength.  Let $F_{e}$ be the area with multiplicities enclosed by  the evolute of $\gamma$. Then
  $$\int_{\gamma}k(s)\,ds-2\pi=c|F_{e}|,$$
  where $k(s)$ is the curvature of $\gamma(s)$ in the ambient space.
  \end{teorema}

The  strong convexity notion used above will be defined later. 

\medskip

As it is well known the isoperimetric inequality in $X_{c}^2$ states
$$F\leq \frac{L^2+cF^2}{4\pi},$$ where $L$ is the length of a simple closed curve $\gamma$ and $F$ the area enclosed by $\gamma$; see for instance \cite{Santa}. We estimate  the 
isoperimetric deficit by means of the area enclosed by the evolute of $\gamma$ proving  the following result. 


\begin{teorema}\label{94} Let $\gamma$ be a positively oriented closed strongly convex  curve on $X^2_c$ of length $L$. Let $F$ be the area enclosed by $\gamma$.   Then the isoperimetric deficit $\Delta=L^2-4\pi F+c F^2$ is bounded by $$ \Delta \leq cF_{e}^2+4\pi |F_{e}|,$$
  where  $F_{e}$ is the area with multiplicities enclosed by the evolute of $\gamma$. Equivalently,
 $$\Delta\leq \frac{1}{c}\biggl(\bigl(\int_{\gamma}k(s)\,ds\bigr)^2-4\pi^2\biggr),$$
where $k(s)$ is the curvature of $\gamma$ in the ambient space. Equality holds if and only if $\gamma$ is a circle.
 \end{teorema}

Finally we provide a Gauss-Bonnet formula with multiplicities (Theorem \ref{GB}) that enables us to calculate the total curvature of the evolute of a curve, for the special case of evolutes with a finite number of singular points, these being the points at which the evolute fails to have a tangent. We prove the following result.

\begin{teorema}\label{150}Let $\gamma$ be a positively oriented closed strongly convex curve on $X^2_c$ and let $\gamma_{e}(s_{e})$ be the evolute of $\gamma$, where $s_{e}$ is its arclength parameter. Assume that $\gamma_{e}(s_{e})$ has a finite number of singular points. Then the integral of the geodesic curvature $k_{e}(s_{e})$  of the  evolute $\gamma_{e}(s_{e})$  is given by 
$$\int_{\gamma_{e}}k_{e}(s_{e})\,ds_{e}=c|F_{e}|+2\pi,$$
 where $F_{e}$ is the area with multiplicities enclosed by $\gamma_{e}$.
\end{teorema}

We point out that the obstruction to generalize the previous result for the evolute of an arbitrary curve comes from the fact that the tangent vector to the evolute can vanish on an arbitrary closed set.  We overcome this difficulty considering only evolutes with a finite number of singularities.


\section{Preliminaries}
We recall here the notions of geodesic curvature and radius of curvature of a curve in $X^2_{c}$.   

In order to treat together the cases of constant positive and negative curvature we consider, as in  \cite{Rat} or \cite{Santa}, the metric on $\R^3$  given by the matrix \begin{eqnarray}\label{a}\left(\begin{array}{ccc}1 & 0 & 0 \\0 & 1 & 0 \\0 & 0 & \epsilon\end{array}\right),\end{eqnarray} where $\epsilon=\pm 1$. If $\epsilon=1$ it is a Riemannian metric and if $\epsilon=-1$ it is a Lorentz metric.

The scalar product of the vectors $u$ and $v$ is denoted by $\langle u,v\rangle.$ The subspace of $\R^3$ given by 
$$S(\epsilon,K)=\{u\in\R^3;\; \langle u,u\rangle =\frac{1}{\epsilon K}\}$$ 
where $K$ is a positive constant, is the standard sphere of radius $R=\frac{1}{\sqrt{K}}$ if $\epsilon=1$ or a hyperboloid if $\epsilon=-1$. In this second case we assume that the elements $u=(u_{1},u_{2},u_{3})$ of $S(\epsilon, K)$ satisfy $u_{3}>0$. Since $S(-1, K)$ consists of vectors of norm $Ri$, it is also called  the imaginary sphere.

In both cases, $\epsilon=1$ or $\epsilon=-1$, $S(\epsilon,K)$ is a Riemannian manifold of constant curvature $c=\epsilon K$. In fact, the metric \eqref{a} restricted to $S(\epsilon, K)$ is positive definite.
Hence $X^2_c=S(1,c)=\Ss^2_c$ for $c>0$, $X^2_c=S(-1,-c)=\hh^2_c$ for $c<0$. 

  We also note that   the tangent space to $S(\epsilon)$ at $P\in S(\epsilon, K)$, $T_{P}S(\epsilon, K)$, is given by $$T_{P}S(\epsilon, K)=P^{\bot}$$
where $P^{\bot}$ denotes the subspace of $\R^3$ orthogonal (with respect to the Riemannian or the Lorentz metric) to $P$.

 Since the covariant derivative on $S(\epsilon, K)$ is the orthonormal projection on $S(\epsilon, K)$ of the covariant derivative of $\R^3$ we have
$$\nabla_{v}Y=v(Y)-c \,\langle v(Y),P\rangle P$$
where $v\in T_{P}S(\epsilon, K)$, $Y$ is a tangent vector field on $S(\epsilon, K)$, and $v(Y)=(v(Y_{1}),v(Y_{2}),v(Y_{3}))$ is the directional derivative of each component. Note that since $\langle \nabla_{v}Y,P\rangle=0,$  we have  $\nabla_{v}Y\in T_{P}S(\epsilon, K)$. 

Let now $\gamma(t)$ be a {\em regular} curve on $S(\epsilon, K)$, that is $\gamma(t)$ is smooth and $\gamma'(t)\neq 0$,  and take $v=Y=\gamma'(t)$.  We have
$$\nabla_{\gamma'(t)}\gamma'(t)=\gamma''(t)-c\, \langle\gamma''(t),\gamma(t)\rangle\gamma(t).$$
If $\gamma$ is parametrized by arclength $s$,  then 
\begin{eqnarray*}
\langle \gamma(s),\gamma(s)\rangle&=&1/c,\qquad \langle \gamma'(s),\gamma(s)\rangle =0,\qquad \langle \gamma'(s),\gamma'(s)\rangle =1,\\
\langle \gamma''(s),\gamma(s)\rangle&=&-1,\qquad
\langle \gamma''(s),\gamma'(s)\rangle =0,
\end{eqnarray*}
and hence
\begin{eqnarray}\label{nabla}\nabla_{\gamma'(s)}\gamma'(s)=\gamma''(s)+ c \,\gamma(s).\end{eqnarray}

\begin{defi}Let $\gamma(s)$ be a regular curve on $X_{c}^2$ parametrized by arclength. The {\em geodesic curvature} $k_{g}(s)$ of $\gamma(s)$  is $$k_{g}(s)=|\nabla_{\gamma'(s)}\gamma'(s)|.$$
The {\em normal vector} $n(s)$ to $\gamma(s)$ 
is given by 
$$\nabla_{\gamma'(s)}\gamma'(s)=k_{g}(s)n(s).$$
\end{defi}

%
%
%
%

Note that, for $c\neq 0$, $n(s)$ is not the principal normal of $\gamma(s)$ as a curve in the ambient space  $\R^3$ (Euclidean or Lorentzian).

We shall use later  the equality 
\begin{eqnarray}\label{gghh}\gamma''(s)=k_{g}(s)n(s)- c \gamma(s).\end{eqnarray}

If the parameter $t$ of a given curve $\gamma(t)$ on $X^2_c$ is not the arclength parameter, 
the geodesic curvature is given by 

\begin{eqnarray}\label{hh}
k_{g}(t)=f(t)^2 \langle \nabla_{\gamma'(t)}\gamma'(t),n(t)\rangle=f(t)^2\langle\gamma''(t),n(t)\rangle,  
\end{eqnarray}
where $f(t)^2=\langle\gamma'(t),\gamma'(t)\rangle^{-1}$ and $n(t)$ is the normal vector to $\gamma(t)$.

The relationship between the geodesic curvature $k_{g}(s)$ and the curvature $k(s)$ of $\gamma(s)$ as a curve in $\R^3$ is 
\begin{eqnarray}\label{curvat}\sqrt{k_{g}^2(s)+c}=k(s),
\end{eqnarray}
since
\begin{eqnarray*}
k_{g}(s)^2=\langle \gamma''(s)+c\gamma(s), \gamma''(s)+c\gamma(s)\rangle=k(s)^2+2c\langle\gamma(s),\gamma''(s)\rangle+c^2\langle\gamma(s),\gamma(s)\rangle=k(s)^2-c.
\end{eqnarray*}

\bigskip

In order to define the radius of curvature, we shall use the generalized sinus and cosinus functions:
\[ \sn\rho := \left\{ \begin{array}{lll}
             \frac{1}{\sqrt{-c}} \sinh (\sqrt{-c}\, \rho) & , & c<0 \\
             \rho &  , & c=0  \\
             \frac{1}{\sqrt{c}} \sin (\sqrt{c}\, \rho) & , & c>0
                           \end{array}  \right.  \]

\[ \cn\rho := \left\{ \begin{array}{lll}
             \cosh (\sqrt{-c}\, \rho) & , & c<0 \\
             1 &  , & c=0  \\
             \cos (\sqrt{c}\, \rho) & , & c>0 \; ,
                           \end{array}  \right.  \]
as well as $\tan_{c}\rho=\dfrac{\sn\rho}{\cn\rho}$ and $\cot_{c}\rho=\dfrac{\cn\rho}{\sn\rho}$.
\bigskip

\begin{defi}
We say that a regular simple curve $\gamma(s)$ on  $X^2_c$  parametrized by arclength is {\em strongly convex} if,  for each $s$,  $k_{g}(s)>0$ for $c\geq 0$ or   $k_{g}(s)>\sqrt{|c|}$ for $c<0$.
\end{defi}

This enable us to give the following definition.
\begin{defi}\label{def23}
Let $\gamma(s)$ be a strongly convex curve on $X^2_c$  parametrized by arclength. The {\em radius of  curvature} of $\gamma(s)$ is the function $\rho(s)$  defined by 
$$k_{g}(s)=\ctn \rho(s),$$
 where  $\ctn \rho(s)$ is the generalized cotangent function.\end{defi}

The condition of strongly convexity corresponds, for $c<0$, to the notion of horocyclic convexity. It is needed because, for $c<0$, $\ctn x>\sqrt{-c}$, for all $x\in\R$. For $c>0$ we shall also assume that $0<\sqrt{c}\rho(s)<\pi/2$. 

The motivation for the Definition \ref{def23} is the fact that a circle of radius $\rho$ has geodesic curvature $\cot_{c}\rho.$

\section{Evolutes}

First we recall that given  $x\in X^2_{c}$ and $y\in T_{x}X^2_c$, with $\langle y,y\rangle=1$, then $$\sigma(t)=\cn(t) \, x+\sn (t)\, y, $$ is the geodesic through $\sigma(0)=x$ with  director tangent vector $\sigma'(0)=y$. This it easy to see, since $\sigma(t) $ verifies the equation of the geodesics 
$\sigma''(t)+c\sigma(t)=0.$ Moreover $t$ is the arclength of $\sigma(t)$ because $\langle\sigma'(t),\sigma'(t)\rangle=1.$

\begin{defi}
Let $\gamma(s)$ be a strongly convex curve on  $X^2_c$ parametrized by arclength.
The {\em evolute} of $\gamma$ is the curve
$$\gamma_{e}(s)=\cn \rho(s)\gamma(s)+\sn\rho(s)n(s)$$
where $\rho(s)$ and $n(s)$  are respectively the radius of curvature and the normal to $\gamma(s)$. 
\end{defi}
So $\gamma_{e}(s)$ is the point on the geodesic through $\gamma(s)$ with director tangent vector $n(s)$,  given by the value $\rho(s)$ of the parameter.  Remark that $s$ is not the arclength parameter of the evolute. 

By the definition of $n(s)$, equation \eqref{nabla},  and the definition of $k_{g}(s)$,   we have
\begin{eqnarray}\label{uu}n(s)=\tan_{c}\rho(s)(\gamma\,''(s)+ c \gamma(s)),\end{eqnarray}
and hence
\begin{eqnarray*}\label{1907}\gamma_{e}(s)=\frac{1}{\cn\rho(s)}\left(\gamma(s)+\sn^2\rho(s)\,\gamma\,''(s)\right).\end{eqnarray*}

For further purposes we need to compute the tangent vector to the evolute.  

We first compute the derivative of the vector $n(s)$.   Since $\langle n(s),n(s)\rangle=1$, we have \linebreak $\langle n\,'(s),n(s)\rangle=0$ and hence
$$n\,'(s)=a(s)\gamma(s)+b(s)\gamma\,'(s).$$
The equality  $\langle \gamma(s),n(s)\rangle=0$ implies $\langle \gamma(s), n\,'(s) \rangle=0$, and so $a(s)=0$. Also, from formula $\eqref{gghh}$, we have
$$k_{g}(s)=\langle\gamma\,''(s),n(s)\rangle=-\langle\gamma\,'(s),n\,'(s)\rangle=-b(s).$$
Thus, \begin{eqnarray}\label{vv}n\,'(s)=-\ctn \rho(s)\,\gamma\,'(s).\end{eqnarray}

The tangent vector to the evolute is given by
\begin{eqnarray}\label{gg}\frac{d\gamma_{e}(s)}{ds}=\gamma_{e}\,'(s)=\rho'(s)\left(-c\sn \rho(s)\gamma(s)+\cn\rho(s)n(s)\right),\end{eqnarray}
because, according to \eqref{vv},  $$\cn\rho(s)\gamma\,'(s)+\sn\rho(s)n\,'(s)=0.$$

Note that, by  \eqref{uu},  \begin{eqnarray}\label{tangent}\gamma_{e}\,'(s)=\rho'(s)\sn \rho(s)\;\gamma\,''(s).\end{eqnarray}

In particular $\gamma_{e}'(s)=0$ at the critical points of $\rho(s)$. Points where $\rho'(s)\neq 0$ are called {\em regular} points of $\gamma_{e}(s)$ and points where $\rho'(s)= 0$ are called {\em singular} points of $\gamma_{e}(s)$. In a neighborhood of each regular point the evolute can be reparametrized by arclength, and so the normal vector is well defined at these points.

We remark that for $c\neq 0$ the tangent vector to the evolute does not coincide with the normal vector to the curve (at corresponding points). 
 Nevertheless we have the following proposition. 
\begin{prop}\label{prop31}The normal vector to the evolute coincides at regular points, up to the sign,  with the tangent vector to the curve at corresponding points. \end{prop}
{\em Proof.} Let  $n_{e}(s)$ be the  normal vector to the evolute at regular points of    $\gamma_{e}(s)$. 

\begin{center}
\includegraphics[width=.7\textwidth]{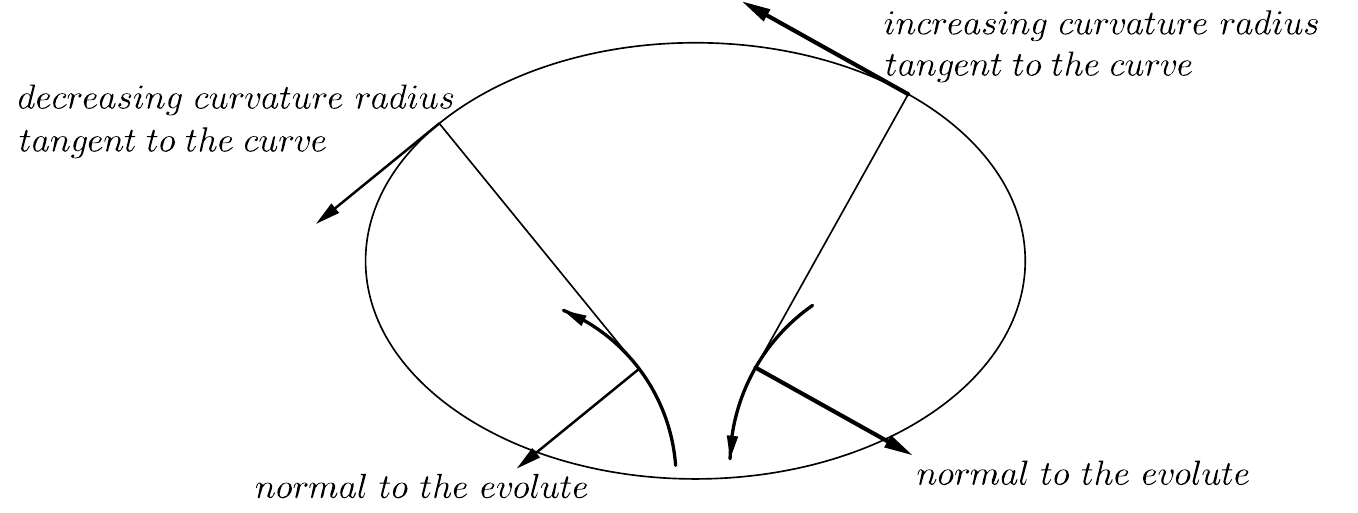}
\end{center}
\centerline{Figure 1.}

We can write  $$n_{e}(s)=A(s)\gamma(s)+B(s)\gamma\,'(s)+C(s)n(s),$$ for some functions $A(s),B(s),C(s)$.
Multiplying by $\gamma_{e}(s)$ one obtains  $$c\, C(s)=-k_{g}(s)A(s)$$ and multiplying by $\gamma_{e}\,'(s)$ one obtains $$A(s)=k_{g}(s)C(s).$$ Since $\gamma(s)$ is strongly convex, we obtain   $A(s)=C(s)=0$, and hence $$n_{e}(s)=B(s)\gamma\,'(s).$$
Thus $|B(s)|=1$ and so $n_{e}(s)=\pm \gamma'(s)$. 

To be more precise, using locally the arclength parameter $s_{e}$ of $\gamma_{e}(s)$, we have  \begin{eqnarray*}B(s)&=&\langle n_{e}(s), \gamma'(s)\rangle=\langle  \frac{d^2\gamma_{e}}{ds_{e}^2},\gamma'(s)\rangle=\langle \frac{d^2\gamma_{e}}{ds^2}(\frac{ds}{ds_{e}})^2, \gamma'(s)\rangle=(\frac{ds}{ds_{e}})^2\langle \rho'(s)\sn\rho(s)\gamma'''(s), \gamma'(s)\rangle\\&=&-(\frac{ds}{ds_{e}})^2\rho'(s)\sn\rho(s)\langle \gamma''(s), \gamma''(s)\rangle.\end{eqnarray*}  Since all the factors  in the right-hand side out of $\rho'(s)$  are positive, we have (see Figure 1)
$$n_{e}(s)=\left\{\begin{array}{l}\hspace{.25cm}\gamma'(s)\quad\mbox{if }\rho'(s)<0\\-\gamma'(s)\quad\mbox{if }\rho'(s)>0.\qquad \quad\square\end{array}\right.$$

We shall need also to compute the geodesic curvature of the evolute of a given curve. Due to equality  \eqref{hh} this notion is well defined at regular points. 
\begin{prop}\label{prop33} The geodesic curvature $k_{e}(s)$ of the evolute of a  strongly convex curve $\gamma(s)$ in $X^2_{c}$,  at regular points, is given by $$k_{e}(s)=\frac{k(s)}{|\rho'(s)|}=\frac{1}{|\rho'(s)|\sn(\rho(s))},$$  
where $k(s)$ is the curvature of $\gamma(s)$ in the ambient space, and $\rho(s)$ is the radius of curvature of $\gamma(s)$. 
\end{prop}
{\em Proof.}
Applying formula \eqref{hh}, Proposition \ref{prop31} and equality \eqref{gg}, we have 
$$k_{e}(s)=\frac{1}{|\gamma_{e}\,'(s)|^2}\langle\gamma_{e}\,''(s),n_{e}(s)\rangle=\pm\frac{1}{\rho'(s)^2}\langle\gamma_{e}\,''(s),\gamma\,'(s)\rangle.$$
Differentiating  the expression of $\gamma_{e}\,'(s)$ obtained in \eqref{gg}, it follows
$$\gamma_{e}\,''(s)=-c\rho'(s)\sn(\rho(s))\gamma\,'(s)+ \rho'(s)\cn(\rho(s))n\,'(s)+\mbox{terms orthogonal to $\gamma\,'(s)$}.$$
Substituting in this expression $n\,'(s)$ by the value obtained in \eqref{vv}, we have $$\gamma_{e}\,''(s)=-\frac{\rho'(s)}{\sn(\rho(s))}\gamma\,'(s)+\mbox{terms orthogonal to $\gamma\,'(s)$}.$$
Hence,

\begin{eqnarray*}k_{e}(s)&=&\pm\frac{1}{\rho'(s)\sn(\rho(s))}.\end{eqnarray*}
Since $k_{e}>0$ we have, 
\begin{eqnarray}\label{kk}k_{e}(s)&=&\frac{1}{|\rho'(s)|\sn(\rho(s))}.\end{eqnarray}
Using the generalized tangent and cotangent functions, 
 it is easy to see
that

 \begin{eqnarray}\label{dd}c \tan_{c}\frac{\rho(s)}{2}=-k_{g}(s)+\sqrt{k_{g}(s)^2+c}\;,\end{eqnarray}
where $k_{g}(s)= \ctn \rho(s)$.
From this and  \eqref{curvat} one obtains
\begin{eqnarray*}\label{ddd}\frac{1}{\sn(\rho(s))}=\sqrt{k_{g}^2(s)+c}=k(s).\end{eqnarray*}
Hence, equation \eqref{kk} can be written as
$$k_{e}(s)=\frac{k(s)}{|\rho'(s)|}.\qquad\square$$


We now  introduce the {\em index} or {\em winding number} of a closed curve on $X_{c}^2$ with respect to a given point. 
%
%
%
%
%
%
%
%

First we recall that the index of a closed piece-wise ${\cal C}^1$ curve $\gamma$ of $\R^2$ is the function defined   by 
$$\mbox{Ind}(\gamma,P)=\frac{\psi_{P}(L)-\psi_{P}(0)}{2\pi},\qquad P\in \R^2\setminus\gamma,$$ where $\psi_{P}(s)$ is a branch of the argument of the vector $(\gamma(s)-P)\in\R^2$, and $s\in [0,L]$ is the arclength parameter of $\gamma$.

It is well known that Ind$(\gamma,P)$ is constant for $P$ in a connected component of $\R^2\setminus \gamma$ and vanishes on the unbounded component. Moreover Ind$(\gamma,P)$ can be computed counting the signed number of intersections of $\gamma$ with a fixed ray starting from $P$; see \cite{BC}, p. 27. 

\medskip

Let now $\gamma$ be a closed curve on $X_{c}^2$ and $P$ a point not on $\gamma$. Assume, without lost of generality, that $\gamma$ and $P$ are contained in an oriented local chart $(U,\varphi)$ where $\varphi:U\longrightarrow X_{c}^2$, and $U$ is an open subset of the plane $\R^2$.  We define Ind$(\gamma, P)$ as Ind$(\tilde{\gamma}, \tilde{P})$, with $\gamma=\varphi\circ\tilde{\gamma}$ and $P=\varphi(\tilde{P})$. 
It is easy to see that this number does not depend on the chosen local chart.   

\begin{defi}
Let $\gamma$ be the a closed piece-wise ${\cal C}^1$ curve on $X_{c}^2$, not necessarily simple. The  {\em area with multiplicities}, $F$, enclosed by $\gamma$ is defined as
 $$F=\int_{X_{c}^2} \operatorname{Ind}(\gamma,P)\;dS,$$
where $dS$ is the area element of $X_{c}^2$.
\end{defi}

\begin{remark}\label{signe}\em Let $\gamma$ be a plane strongly convex closed curve,  positively oriented.  This means Ind$(\gamma,P)=1$ for $P$ in the interior of $\gamma$. Let $\gamma_{e}$ denote the evolute of $\gamma$ and $F_{e}$ the area with multiplicities enclosed by $\gamma_{e}$. We shall see that $F_{e}\leq 0$, a fact that comes from the inequality  \begin{eqnarray}\label{www}\mbox{Ind}(\gamma,P)\cdot \mbox{Ind}(\gamma_{e},P)\leq 0.\end{eqnarray}
Indeed, if $P$ does not belong to a bounded component of $\R^2\setminus\gamma$ or of $\R^2\setminus\gamma_{e}$ at least one of the two indices are zero and the inequality holds.
On the other case, for a fixed $s$,   we have
\begin{eqnarray*}
\gamma(s)-P&=&a\gamma'(s)+bn(s), \quad b\leq 0,\\
\gamma_{e}(s)-P&=&c\gamma_{e}'(s)+dn_{e}(s), \quad d\geq 0,
\end{eqnarray*}
and by Proposition \ref{prop31} 
$$d=\langle \gamma_{e}(s)-P, n_{e}(s)\rangle=\langle \gamma(s)+\rho(s)n(s)-P, n_{e}(s)\rangle= a\langle \gamma'(s), \pm \gamma'(s)\rangle.$$ 

More precisely, 
\begin{eqnarray*}
d&=&a \;\;\,\quad\mbox{  if  }\rho'(s)< 0,\\
d&=&-a \quad\mbox{  if  }\rho'(s)> 0.
\end{eqnarray*}
Since $d> 0$, we have $a\rho'(s)< 0$.

It follows easily that 
\begin{eqnarray*}
\det(\gamma_{e}(s)-P,\gamma_{e}'(s))=a\rho'(s)\det(\gamma'(s), n(s))=a\rho'(s)<0,
\end{eqnarray*} and the inequality \eqref{www}  is proved.  
Note that $\det(\gamma'(s),n(s))=1$ because $\gamma$ is positively oriented.  

From this and the definition of the index of the evolute of a closed strongly convex curve in $X_{c}^2$ it follows that Ind$(\gamma_{e},P)\leq 0$, for $P\in\R^2\setminus\gamma_{e}$. So the area with multiplicities, $F_{e}$, enclosed by $\gamma_{e}$ is negative or zero. \end{remark}

\section{Area of the evolute and total curvature}
  We begin with some notation and a technical lemma.
  Let $\gamma(s)$ be a strongly convex curve on $X^2_c$  parametrized by arclength $s$. 

At each point $\gamma_{e}(s)$ of the evolute  of $\gamma(s)$ we consider the  vector  $T(s)\in T_{\gamma_{e}(s)}X^2_c$ given by $$T(s)=c\,\sn \rho(s)\,\gamma(s)-\cn\rho(s)\,n(s),$$
where $n(s)$ is the normal vector to  $\gamma(s)$.

Note that $T(s)$ is a vector field along  $\gamma_{e}(s)$  which  by \eqref{gg} has  the same direction  than the tangent vector to the evolute at regular points,  but with the advantage that it is also defined at singular points.

We denote, as usual, $$\frac{DT(s)}{ds}\in T_{\gamma_{e}(s)}X^2_c$$ the covariant derivative of $T(s)$ along  $\gamma_{e}(s)$.  For $c\neq 0$, it is the projection on $S(1,c)$ or $S(-1,-c)$ of the directional derivative on $\R^3$.

\begin{lema}\label{lema111} Let $\gamma(s)$ be a strongly convex  curve on $X^2_c$  parametrized by arclength $s$. Then 
$$k(s)=\langle \frac{DT(s)}{ds},\gamma\,'(s)\rangle, $$
where $k(s)$ is the curvature of $\gamma(s)$ in the ambient space.
\end{lema}
{\em Proof.} Since $$\gamma\,'(s)\in T_{\gamma_{e}(s)}X^2_c,$$ we have\begin{eqnarray*}
\langle \frac{DT(s)}{ds},\gamma\,'(s)\rangle&=&\langle \frac{dT(s)}{ds},\gamma\,'(s)\rangle\\&=&\langle c\rho'(s)(\cn\rho(s)\gamma(s)+\sn\rho(s)n(s)),\gamma\,'(s)\rangle\\&+&\langle c\sn\rho(s)\gamma\,'(s)-\cn\rho(s)n\,'(s),\gamma\,'(s)\rangle.
\end{eqnarray*}

By equation \eqref{vv} and Proposition \ref{prop33} we have\begin{eqnarray*}
\langle \frac{DT(s)}{ds},\gamma\,'(s)\rangle&=&c\sn\rho(s)+\cn\rho(s)\ctn\rho(s)=\frac{1}{\sn\rho(s)}=k(s),
\end{eqnarray*}
and lemma is proved. $\square$

\bigskip
 Next result can be seen as a sort of refinement  of Fenchel's Theorem.  
  
  \begin{teorema}\label{teo61}
  Let $\gamma(s)$ be a positively oriented  closed strongly convex  curve on $X^2_c$ parametrized by arclength.  Let $F_{e}$ be the area with multiplicities of the evolute of $\gamma$. Then
  $$\int_{\gamma}k(s)\,ds-2\pi=c|F_{e}|,$$
  where $k(s)$ is the curvature of $\gamma(s)$ in the ambient space.
  \end{teorema}
{\em Proof.} Let $(e_{1},e_{2})$ be a local orthonormal frame of vector fields  on $ X^2_c$. The connection $1$-form $\omega_{12}$ associated to this moving frame is given by
$$\omega_{12}(X)=\langle\nabla_{X}e_{1},e_{2}\rangle$$
for each tangent vector field $X$.

In the vector tangent space $T_{\gamma_{e}(s)}X^2_c$ we have

$$T(s)=\cos\theta(s)e_{1}+\sin\theta(s)e_{2}$$
where $\theta(s)$ is the angle, module $2\pi$, between $T(s)$ and $e_{1}$.

Then, by Lemma \ref{lema111}, 
\begin{eqnarray*}k(s)&=&\langle \frac{DT(s)}{ds},\gamma\,'(s)\rangle=\langle \frac{D(\cos\theta(s)e_{1}+\sin\theta(s)e_{2})}{ds},\gamma\,'(s)\rangle\\&=&\langle \theta'(s)(-\sin\theta(s)e_{1}+\cos\theta(s)e_{2}),\gamma\,'(s)\rangle+\langle \cos\theta(s)\frac{De_{1}}{ds}+\sin\theta(s)\frac{De_{2}}{ds},\gamma\,'(s)\rangle.\end{eqnarray*}
But $$\gamma\,'(s)=-\sin\theta(s)e_{1}+\cos\theta(s)e_{2}\in T_{\gamma_{e}(s)}X^2_c$$
and
\begin{eqnarray*}
\langle\frac{De_{1}}{ds}, -\sin\theta(s)e_{1}+\cos\theta(s)e_{2}\rangle&=&\cos\theta(s)\omega_{12}(\gamma_{e}\,'(s))
\\
\langle\frac{De_{2}}{ds}, -\sin\theta(s)e_{1}+\cos\theta(s)e_{2}\rangle&=&\sin\theta(s)\omega_{12}(\gamma_{e}\,'(s)). 
\end{eqnarray*}
Hence
\begin{eqnarray*}k(s)&=&\theta'(s)+\omega_{12}(\gamma_{e}\,'(s)).\end{eqnarray*}
This yields to an equality of $1$-forms
$$k(s)ds=d\theta+\gamma_{e}^*\;\omega_{12}$$
Integrating on $[0,L]$ we have,

$$\int_{[0,L]}k(s)ds=\int_{[0,L]} d\theta+\int_{[0,L]} \gamma_{e}^*\;\omega_{12}.$$
Equivalently, 
$$\int_{0}^L k(s)\,ds=\int_{0}^L \theta'(s)\,ds+\int_{\gamma_{e}} \omega_{12}$$
But we know, from the structure equations (see, for instance, \cite{SPI1970}, Vol. II, p. 295), that $$d\omega_{12}=-c\,\theta^1\wedge\theta^2=-c\, dS$$
where $(\theta^1,\theta^2)$ is the dual basis of $(e_{1},e_{2})$ and $dS$ the area element of $X^2_c$. 

By the  Green formula with multiplicities (see for instance \cite{BC}, p. 213) we have
\begin{eqnarray*}\int_{0}^L k(s)\,ds&=&\int_{0}^L \theta'(s)\,ds+\int_{X^2_c} \mbox{Ind}(\gamma_{e},P)\,d\omega_{12}\\&=&2\pi-c\int_{X^2_c} \mbox{Ind}(\gamma_{e},P)dS\\&=&2\pi+c|F_{e}|\end{eqnarray*}
since the index of the evolute is negative (see remark \ref{signe}), and theorem is proved. $\square$

\bigskip

Next we give, using Theorem \ref{teo61}, a simple proof of a known result which appears in \cite{ERS} (Theorem 3.8) but with a completely different proof. It will be used in Section 5.

\begin{teorema}\label{teopral}Let $\gamma(s)$ be a positively oriented closed strongly convex  curve on $X^2_c$  parametrized by arclength. Let  $\rho(s)$  be the corresponding radius of curvature.  Then
$$\int_{\gamma}\tan_{c}\frac{\rho(s)}{2}ds=F+|F_{e}|,$$
where $F$ is the area enclosed by $\gamma$ and $F_{e}$ is the area with  multiplicities enclosed by the evolute of $\gamma$.
\end{teorema}
{\em Proof.}
%
%
%
Integrating both sides of \eqref{dd} and using \eqref{curvat} one obtains
$$c\int_{\gamma}\tan_{c}\frac{\rho(s)}{2}ds=-\int_{\gamma} k_{g}(s)ds+\int_{\gamma} k(s)\;ds.$$
By the Gauss-Bonnet theorem  (see for instance \cite{Santa}, p. 303) and Theorem \ref{teo61} we have
 
$$c\int_{\gamma}\tan_{c}\frac{\rho(s)}{2}ds=(-2\pi+c F)+(2\pi+c |F_{e}|)=c F+c |F_{e}|, $$
and theorem is proved. $\square$

As an immediate consequence we have the following Corollary, that can be considered as a generalization to the case of constant curvature of the  $2$-dimensional analogue of Ros' inequality; see \cite{ERev}.

\begin{coro}\label{coro62}Let $\gamma(s)$ be a positively oriented closed  strongly convex  curve on $X^2_c$  parametrized by arclength. Let $\rho(s)$  be the radius of curvature of $\gamma(s)$. Then
$$F\leq \int_{\gamma}\tan_{c}\frac{\rho(s)}{2}ds,$$
where $F$ is the area enclosed by $\gamma$. Equality holds if and only if $\gamma$ is a circle.
\end{coro}
{\em Proof.} The inequality is immediate from Remark \ref{signe} and Theorem \ref{teopral}. Equality holds if and only if $F_{e}=0$. Since Ind$(\gamma_{e},P)\leq 0$ (see remark \ref{signe}), it must be Ind$(\gamma_{e},P)=0$. This implies that the evolute $\gamma_{e}$ is a point and hence $\gamma$ must be a  circle. Indeed, if the evolute $\gamma_{e}$ was not  a point we could choose a small ball separated by $\gamma_{e}$ in two connected components. Then the index would be a different integer in each of these parts since although the evolute can  be traversed twice this always happens in the same sense. This gives a contradiction. $\square$


\bigskip

Since the evolute of a simple closed curve $\gamma$ coincides with the evolute of a curve `parallel' to it, the above results relating the area enclosed by $\gamma$ and the area enclosed by its evolute  yield a new proof of Steiner's formula for tubes on noneuclidean spaces; see  \cite{Santa}, p. 322.

\begin{teorema}[Steiner formula]\label{Steiner} Let $\gamma=\partial Q$ be the strongly convex boundary of a  compact domain $Q$ in  $X^2_{c}$. Denote by $F$ the area of $Q$ and by $L$ the length of $\gamma$. 
Let $Q_{r}$ be the semitube around $Q$ in the direction of the outward normal.
 Then \begin{eqnarray*}F_{r}-F=L\sn (r) +2\sn^2(r/2)(2\pi-cF)\end{eqnarray*} where $F_{r}$ denotes the area of $Q\cup Q_{r}$.
\end{teorema}
{\em Proof.}
Applying Theorem \ref{teopral} to $\gamma$ and to $\gamma_{r}=\partial(Q\cup Q_{r})$, and taking into account that the evolute of $\gamma$ coincides with the evolute of $\gamma_{r}$, and that the curvature radius of $\gamma_{r}$ and $\gamma$, at corresponding points $\gamma(s)$ and $\gamma_{r}(s)=\mbox{exp}_{\gamma(s)}rN(s)$, are related by $\rho_{r}(s)=\rho(s)+r$, 
 we have
 
$$F_{r}-F=\int_{\gamma_{r}}\tan_{c}\frac{\rho(\tau)+r}{2}d\tau-\int_{\gamma}\tan_{c}\frac{\rho(s)}{2}ds$$ 
 where $ds$ is the arclength measure on $\gamma$, and $d\tau$ is the arclength measure on $\gamma_{r}$.

Applying the sinus theorem in the infinitesimal triangle of the Figure 2
we see that
 $$d\tau=\frac{\sn(\rho(s)+r)}{\sn \rho(s)}ds.$$

\pagebreak
\begin{center}
\includegraphics[width=.7\textwidth]{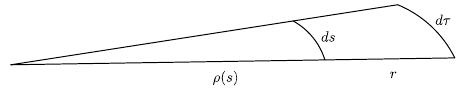}
\end{center}
\centerline{Figure 2.}

 Hence
 
$$F_{r}-F=\int_{\gamma}\left(\tan_{c}\dfrac{\rho(s)+r}{2}\cdot\dfrac{\sn\dfrac{\rho(s)+r}{2}\cn\dfrac{\rho(s)+r}{2}}{\sn\dfrac{\rho(s)}{2}\cn\dfrac{\rho(s)}{2}}-\dfrac{\sn\dfrac{\rho(s)}{2}}{\cn\dfrac{\rho(s)}{2}}\right)ds.$$ 
Simplifying
$$F_{r}-F=\int_{\gamma}\left(\dfrac{\sn^2((\rho(s)+r)/2)-\sn^2((\rho(s)+r)/2)}{\sn(\rho(s)/2)\cn(\rho(s)/2)}\right)ds.$$

Now we substitute  $\sn^2((\rho(s)+r)/2)$ for his expression 
\begin{eqnarray*}\sn^2((\rho(s)+r)/2)&=&\sn^2(\rho(s)/2)\cn^2(r/2)+\cn^2(\rho(s)/2)\sn^2(r/2)\\&+&2\sn(r/2)\cn(r/2)\sn(\rho(s)/2)\sn(\rho(s)/2)\end{eqnarray*}
and we obtain

\begin{eqnarray*}F_{r}-F&=&L\sn (r) +2\sn^2(r/2)\int_{\gamma}\frac{\cn(\rho(s))}{\sn(\rho(s))}ds\\&=&L\sn (r) +2\sn^2(r/2)(2\pi-cF). \quad\square\end{eqnarray*}

\section{An estimate of the isoperimetric deficit}
As it is well known the isoperimetric inequality in $X_{c}^2$ states that
$$F\leq \frac{L^2+cF^2}{4\pi},$$ where $L$ is the length of a simple closed curve $\gamma$ and $F$ the area enclosed by $\gamma$. 

Here we apply previous results to provide an upper bound for the right-hand side of this inequality. 
\begin{teorema}\label{teo71}
Let $\gamma(s)$ be a positively oriented  closed strongly convex  curve on $X^2_c$ of length $L$ parametrized by arclength. Let  $\rho(s)$  be the corresponding radius of curvature.  Then
\begin{eqnarray}\label{kl3}\frac{L^2+cF^2}{4\pi}\leq \int_{\gamma}\tan_{c}\frac{\rho(s)}{2}ds+\frac{cF_{e}^2}{4\pi},\end{eqnarray} where $F$ is the area enclosed by $\gamma$ and  $F_{e}$ is the area with multiplicities enclosed by  the evolute of $\gamma$.
Equality holds if and only if $\gamma$ is a circle. \end{teorema} 
{\em Proof.} Integrating both sides of the identity
$$\ctn \frac{\rho(s)}{2}=\ctn\rho(s)+\frac{1}{\sn \rho(s)}$$
and multipliying by $$\int_{\gamma}\tan_{c}\frac{\rho(s)}{2}ds, $$ we obtain 
$$\int_{\gamma}\tan_{c}\frac{\rho(s)}{2}ds\cdot \int_{\gamma}\ctn\frac{\rho(s)}{2}ds=\int_{\gamma}\tan_{c}\frac{\rho(s)}{2}ds\left(\int_{\gamma} \ctn\rho(s)ds+  \int_{\gamma} \frac{1}{\sn \rho(s)}ds\right)$$

On the other hand, by the Schwarz's inequality,  we have
$$L^2=\left(\int_{\gamma} \sqrt{\tan_{c}\frac{\rho(s)}{2}}\frac{1}{\sqrt{\tan_{c}\dfrac{\rho(s)}{2}}}ds\right)^2\leq \int_{\gamma}\tan_{c}\frac{\rho(s)}{2}ds\cdot \int_{\gamma}\ctn\frac{\rho(s)}{2}ds.$$
Hence, using the Gauss-Bonnet theorem, and Theorems \ref{teo61} and \ref{teopral}, we obtain 

$$L^2\leq (F+|F_{e}|)\bigl((2\pi-c F)+(2\pi+c|F_{e}|)\bigr)=(F+|F_{e}|)\bigl(4\pi-c(F-|F_{e}|)\bigr).$$

Thus
\begin{eqnarray*}\label{kl}L^2\leq 4\pi \int_{\gamma}\tan_{c}\frac{\rho(s)}{2}ds-c(F^2-F_{e}^2),\end{eqnarray*}
and inequality \eqref{kl3} is proved.

Finally, note that equality holds if and only if $k_{g}$ is constant. But closed curves on $X^2_c$ of constant geodesic curvature are circles. $\square$

\medskip
As a consequence  we have an estimate of the isoperimetric deficit in terms of $F_{e}$.
\begin{teorema}\label{94} Let $\gamma$ be a positively oriented  closed strongly convex  curve on $X^2_c$ of length $L$. Let $F$ be the area enclosed by $\gamma$.   Then the isoperimetric deficit $\Delta=L^2-4\pi F+c F^2$ is bounded by $$ \Delta \leq cF_{e}^2+4\pi |F_{e}|,$$
  where  $F_{e}$ is the area with multiplicities enclosed by the evolute of $\gamma$. Equivalently,
 $$\Delta\leq \frac{1}{c}\biggl(\bigl(\int_{\gamma}k(s)\,ds\bigr)^2-4\pi^2\biggr),$$
where $k(s)$ is the curvature of $\gamma$ in the ambient space. Equality holds if and only if $\gamma$ is a circle.
 \end{teorema}
{\em Proof.} First inequality follows from Theorem \ref{teo71} and Theorem  \ref{teopral}, and for the second one we use Theorem \ref{teo61}. $\square$
%
%

\begin{remark}\em 
Combining the isoperimetric inequality and 
formula \eqref{kl3} one gets

$$F\leq \int_{\gamma}\tan_{c}\frac{\rho(s)}{2}ds+\frac{cF_{e}^2}{4\pi},$$
 which is,  for the case $c< 0$,  an improvement of Corollary \ref{coro62}. $\square$
\end{remark}
 
\section{The Gauss-Bonnet theorem for evolutes}

It is possible to have a regular curve with an arbitrary closed set (for instance, a Cantor set) of maximums 
and minimums of its curvature. In this case its evolute has a singular point  corresponding to each point of this closed set. The angle 
between the tangent vector to the evolute and a given direction is not well defined at singular points, since at these points the tangent vector to the evolute vanishes. This is  an obstruction in order to find a formula  for the integral of the geodesic curvature of the evolute. Nevertheless we think that it is interesting to consider the particular case of 
evolutes with a finite number of singular points.

More generally, let us consider a closed piece-wise ${\cal C}^2$ curve $\gamma(s)$ on $X^2_c$ where $s$ is the arclength parameter. That is, $\gamma(s)$ has two  continuous derivatives  except (possibly) at a finite number of singular points at which left and right derivatives exist. The geodesic curvature of $\gamma(s)$ is defined out of these singular points.

For this class of curves we give an extension  of the Gauss-Bonnet theorem.

\begin{teorema}[Gauss-Bonnet theorem with multiplicities]\label{GB}
Let $\gamma(s)$ be a positively oriented closed piece-wise ${\cal C}^2$ curve on $X_{c}^2$, not necessarily simple, where $s$ is the arclength parameter. Then the integral of the geodesic curvature $k_{g}(s)$  is given by $$\int_{\gamma}k_{g}(s)\,ds=-c F+\sum_{k=1}^N\theta_{k}+(2\nu-N)\pi,$$
where $F$ is the area with  multiplicities enclosed by $\gamma$, $N$ is the number of singular points,  $\theta_{k}$  are the interior angles at these points  and $\nu\in \Z$.  

\end{teorema}
{\em Proof.} Suppose that $(u,v)$ is a system of orthogonal coordinates defined on $X^2_c$ given by a parametrization\label{para} $\varphi:U\longrightarrow X^2_c$ defined on an open subset $U$ of the $(u,v)$ plane $\R^2$.  We may assume $\gamma(s)\subset\varphi(U)$ for all $s\in[0,L]$.

If we write the metric in this coordinates as $$\left(\begin{array}{cc}E & 0 \\0 & G\end{array}\right),$$  
the geodesic curvature of the curve $\gamma(s)=\varphi(u(s),v(s))$ is given by the piece-wise  ${\cal C}^1$ function 
$$k_{g}(s)=\frac{1}{2\sqrt{EG}}\left(G_{u}\frac{dv}{ds}-E_{v}\frac{du}{ds}\right)+\frac{d\theta}{ds}$$ where $\theta(s)$ is the the positive angle between $\frac{\partial}{\partial u}_{|\gamma(s)}$ and $\gamma'(s)$.

If we consider the $1$-form on $U\subset\R^2$,  $\omega=Adu+Bdv$, with
$$A=-\frac{E_{v}}{2\sqrt{EG}}, \qquad B=\frac{G_{u}}{2\sqrt{EG}}$$
we have the equality of $1$-forms \begin{eqnarray}\label{1formes}k_{g}(s)ds=\omega+d\theta.\end{eqnarray}
In this equality it is assumed that  $\omega$ is restricted to $\gamma(s)$, and  $d\theta=\frac{d\theta}{ds}ds=\theta'(s) ds$.

On the other hand, it is known that the   Gauss curvature $c$ of $X^2_c$ is given by
$$c=-\frac{1}{2\sqrt{EG}}\left(\left(\frac{E_{v}}{\sqrt{EG}}\right)_{v}+\left(\frac{G_{u}}{\sqrt{EG}}\right)_{u}\right),$$
and hence 
\begin{eqnarray*}d\omega&=&-(\frac{\partial A}{\partial v}-\frac{\partial B}{\partial u})\, du\wedge dv=\left(\left(\frac{E_{v}}{2\sqrt{EG}}\right)_{v}+\left(\frac{G_{u}}{2\sqrt{EG}}\right)_{u}\right)du\wedge dv\\&=&-c\sqrt{EG}du\wedge dv=-c dS.\end{eqnarray*}

The Green formula with multiplicities (see for instance \cite{BC}, p. 235) states 
\begin{eqnarray*}\label{green}\int_{\gamma}\omega=\int_{\R^2}\mbox{Ind}(\gamma, P)\; d\omega\end{eqnarray*}
where Ind$(\gamma, P)$ denotes the  index of the curve $\gamma(s)=\varphi^{-1}(\gamma(s))$ with respect to the point $P$, and $\omega$ is  a $1$-form on $\R^2$.

Hence, integrating both sides of \eqref{1formes},  we have, 

$$\int_{\gamma}k_{g}(s)\,ds=\int_{\gamma}\omega +\int_{\gamma}d\theta=\int_{\R^2}\mbox{Ind}(\gamma,P)\;d\omega+\int_{\gamma}d\theta=-c\int_{\R^2}\mbox{Ind}(\gamma, P)\;dS+\int_{\gamma}d\theta,$$

and since by definition $$F=\int_{\R^2}\mbox{Ind}(\gamma, P)\;dS,$$
we have
\begin{eqnarray}\label{ppp}\int_{\gamma}k_{g}(s)\,ds=-c F+\int_{\gamma}d\theta.\end{eqnarray}

But
\begin{eqnarray*}
\int_{\gamma}d\theta=\sum_{k=0}^{N}\int_{a_{k}}^{a_{k+1}}\theta'(s)ds
\end{eqnarray*}
with $0=a_{0}<a_{1}<\dots<a_{N}<a_{N+1}=L$, where $a_{1},a_{2},\dots, a_{N}$ are the singular points of $\gamma(s)$ and $\gamma(a_{0})=\gamma(L)$,  ($L$ the length of $\gamma$). 
Hence (see Figure 3)
\begin{eqnarray*}\int_{\gamma}d\theta&=&\sum_{k=0}^{N}(\theta(a_{k+1}^-)-\theta(a_{k}^+))\\&=&\sum_{k=1}^{N}(\theta(a_{k}^-)-\theta(a_{k}^+))+(\theta(a_{0}^-)-\theta(a_{N+1}^+))\\&=&-\sum_{k=1}^{N}(\pi-\theta_{k})+2\pi\nu,\qquad \nu\in\Z,\end{eqnarray*}
since by definition of interior angle 
\begin{eqnarray}\label{tangent2}
\theta(a_{k}^+)-\theta(a_{k}^-)=\pi-\theta_{k}. 
\end{eqnarray}
\begin{center}
\includegraphics[width=.6\textwidth]{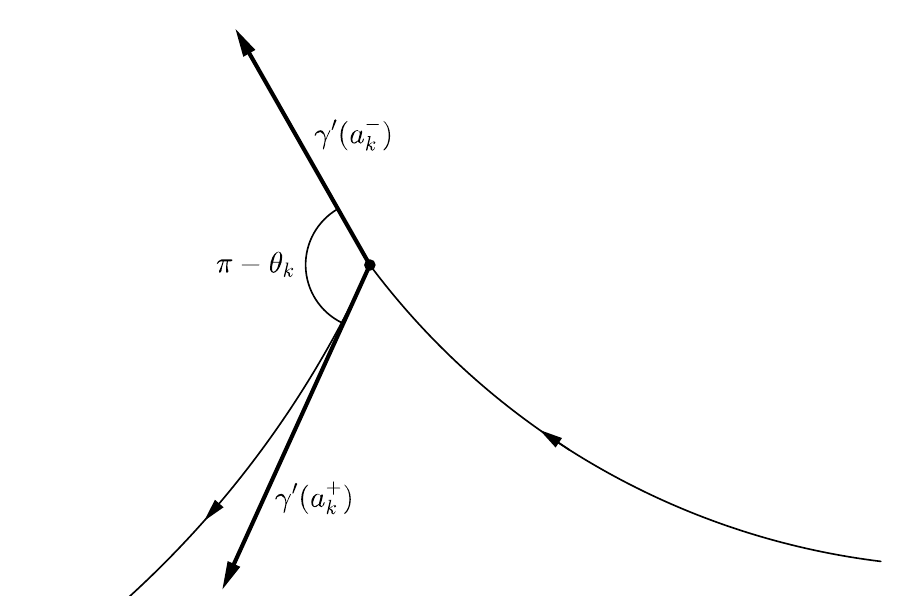}
\end{center}
\centerline{Figure 3.}
\medskip
Substituting this expression of $\int_{\gamma}d\theta$ in \eqref{ppp} the  theorem is proved. $\square$

\medskip

Note that for a plane curve the integer number $\nu$ coincides with its rotation index. Recall that the rotation index of a closed plane curve is defined as the number of turns made by the tangent vector to this curve; see a precise definition in \cite{Chern}. 
\medskip

Using the previous theorem we can compute now the total geodesic curvature of the evolute $\gamma_{e}$ of a strongly convex curve on $X_{c}^2$ in the case that $\gamma_{e}$ has  a finite number of singular points, obtaining a Gauss-Bonnet formula for these evolutes.  Indeed,  we can reparametrize $\gamma_{e}$ with respect to its arclength parameter $s_{e}$ obtaining a piece-wise ${\cal C}^2$ curve to which  Theorem \ref{GB} can be applied. It does not seem possible to do this in the general case.

\begin{teorema}\label{15}Let $\gamma$ be a positively oriented closed strongly convex curve on $X^2_c$ and assume that its evolute $\gamma_{e}$ has a finite number of singular points. Let $s_{e}$ be the arclength parameter of $\gamma_{e}$.  Then the integral of the geodesic curvature $k_{e}(s_{e})$  of the  evolute $\gamma_{e}(s_{e})$  is given by 
$$\int_{\gamma_{e}}k_{e}(s_{e})\,ds_{e}=c|F_{e}|+2\pi,$$
 where $F_{e}$ is the area with multiplicities enclosed by $\gamma_{e}$.
\end{teorema}
{\em Proof.} For a negatively  oriented closed piece-wise ${\cal C}^2$  curve we have, by Theorem \ref{GB},
$$\int_{\gamma}k_{g}(s)\,ds=-c F-\sum_{k=1}^N\theta_{k}+(N+2\nu)\pi.$$

This equality can be applied to $\gamma_{e}(s_{e})$ which is piece-wise ${\cal C}^2$ and negatively oriented by Remark \ref{signe}. To evaluate the right-hand side of previous equality, when applied to $\gamma_{e}$, we consider first of all the case of plane cuves. 

Note that the interior angles $\theta_{k}$ are zero.  This is a consequence of equalities \eqref{tangent} and \eqref{tangent2} and the fact that the angles $\theta(a_{k}^+)$ and $\theta(a_{k}^-)$ in \eqref{tangent2} are the angles with respect to a given direction of the normal vector to the curve and its opposite, respectively.

Applying the turning tangents theorem  to the evolute, see for instance \cite{Chern},  one has
$$2\pi\nu=V_{e}-N\pi,$$
where $V_{e}$ is the differentiable variation of the angle formed by the tangent to the evolute with a given direction (sum of the  variations in each interval where the evolute is regular) and $N$ is the number of critical points of the radius of curvature of $\gamma$. Since the tangent to the evolute coincides up to the sign  with the normal  to the curve we get $V_{e}=2\pi$. Hence $N+2\nu=2$ and the thorem is proved for $c=0$.

To generalize the above arguments to the case  $c\neq 0$ we can argue as follows.    
Let $\varphi_{t}:X_{c}^2\longrightarrow X_{(1-t)c}^2$, for $0\leq t \leq 1$, be  a continuous family of  mappings, $\varphi_{0}$ being the identity  and $\varphi_{1}$ the stereographic projection. For each $t$  consider $\varphi_{t}(\gamma)$ and its corresponding evolute (which is not $\varphi_{t}(\gamma_{e})$).  Since the rotation index $\nu$ of this family of evolutes depends continuously on $t$ and takes integer values, it must be constant. So $N+2\nu=2$ holds,  and the proof is finished.   $\square$

\bigskip

{\em Acknowledgements}. The authors are grateful to Gil Solanes for many helpful conversations during the preparation of this work.

\bibliographystyle{plain}
\bibliography{BibliotecaGNE}

\noindent {\em Departament de Matemàtiques \\Universitat Aut\`{o}noma de Barcelona\\ 08193 Bellaterra, Barcelona\\Catalonia\\

\noindent jcufi@mat.uab.cat, agusti@mat.uab.cat.

}

 \end{document}